\documentclass[11pt]{article}
 \usepackage{cite}
 \usepackage{mathrsfs}
 \usepackage{amsfonts}
 \usepackage{amsmath}
 \usepackage{amsfonts,amssymb,color}
 \usepackage{dsfont}
 \usepackage{curves}
 \usepackage{mathrsfs}
 \usepackage{pifont}
 \usepackage{enumitem}
 \usepackage{amssymb}
 \usepackage{latexsym,amsmath,amssymb,amsfonts,epsfig,graphicx,cite,psfrag}
 \usepackage{eepic,color,colordvi,amscd}
 \usepackage[colorlinks=true,
 linkcolor=blue,
 citecolor=blue,
 urlcolor=blue]{hyperref}

 \newtheorem{theorem}{Theorem}

 \newtheorem{lemma}{Lemma}
 \newtheorem{corollary}{Corollary}

 \newtheorem{claim}{Claim}
 \newtheorem{conjecture}{Conjecture}

 \newcommand{\qed}{\hfill\rule{0.5em}{0.809em}}

 \def\emptyset{\mbox{{\rm \O}}}

 \def\qed{\hfill \ensuremath{\square}}
 
 \def\pf{\noindent {\it Proof. }}
 
\begin{document}
 	
 	\title{Structure, Coloring, and Perfect Divisibility of $(P_2\cup P_4, C_3)$-Free Graphs\thanks{ Partially supported by the NSFC under Grant No.~12601670.}}
 	\author{
 		Di Wu$^{1,}$\footnote{Email: diwu@njit.edu.cn.}  \; Xiaowen Zhang$^{2,}$\footnote{Email: xiaowzhang0128@126.com.}\\
 		\small $^1$Department of Mathematics and Physics\\
 		\small Nanjing Institute of Technology, Nanjing 211167, Jiangsu, China\\
 		\small $^2$Department of Mathematics\\
 		\small East China Normal University, Shanghai, 200241, China
 	}
 	\date{}
 	\maketitle

 	\begin{abstract}
 		
 		Goedgebeur and Schaudt [J. Graph Theory 87 (2018), 188-207] conjectured that every $4$-vertex-critical $(P_7,C_3)$-free graph belongs to a family of seven explicitly defined graphs. In this paper, we establish a structural theorem for connected $(P_2\cup P_4,C_3)$-free graphs. As a consequence, we prove that the Mycielski-Grötzsch graph is the unique $4$-vertex-critical graph in this class, thereby confirming the conjecture of Goedgebeur and Schaudt for $(P_2\cup P_4,C_3)$-free graphs. Our structural theorem also yields a characterization of the chromatic number of these graphs and an $O(n^4)$-time algorithm for deciding whether an $n$-vertex $(P_2\cup P_4,C_3)$-free graph is $3$-colorable.
 		
 		We further study perfect divisibility in the larger class of $(P_2\cup P_4,\text{bull})$-free graphs. We prove that a $(P_2\cup P_4,\text{bull})$-free graph is perfectly divisible if and only if it is Mycielski-Grötzsch graph-free. This result generalizes the main theorem of Deng and Chang [Graphs Combin. 41 (2025), 63].
 	\end{abstract}

 	\begin{flushleft}
 	\end{flushleft}
 	
 	\section{Introduction}\label{introdction}

 	In this paper, all graphs are finite and simple. Let $P_k$ and $C_k$ denote a {\em path} and a {\em cycle} on $k$ vertices, respectively. We say that a graph $G$ {\em contains} a graph $H$ if $G$ has an induced subgraph isomorphic to $H$. A graph $G$ is $H$-{\em free} if it does not contain $H$. More generally, for a family $\mathcal{H}$ of graphs, we say that $G$ is $\mathcal{H}$-{\em free} if $G$ contains no induced subgraph isomorphic to any member of $\mathcal{H}$. 
 	
 	For two vertex-disjoint graphs $G_1$ and $G_2$, their {\em union}, denoted by $G_1\cup G_2$, is the graph with vertex set
 	$
 	V(G_1\cup G_2)=V(G_1)\cup V(G_2)
 	$
 	and edge set
 	$
 	E(G_1\cup G_2)=E(G_1)\cup E(G_2).
 	$
 	
 	For a subset $A$ of $V(G)$ and a vertex $u\in V(G)\setminus A$, we say that $u$ is \emph{complete} to $A$ if $u$ is adjacent to every vertex of $A$, and \emph{anticomplete} to $A$ if $u$ is adjacent to no vertex of $A$. For two disjoint subsets $A$ and $B$ of $V(G)$, we say that $A$ is \emph{complete} to $B$ if every vertex of $A$ is complete to $B$, and \emph{anticomplete} to $B$ if every vertex of $A$ is anticomplete to $B$.
 	
 	For $X\subseteq V(G)$, let $G[X]$ denote the subgraph of $G$ induced by $X$. For a vertex $v\in V(G)$ and a subset $X\subseteq V(G)$, let $N_G(v)$ denote the set of neighbors of $v$ in $G$. Let $d_G(v)=|N_G(v)|$, $M_G(v)=V(G)\setminus\bigl(N_G(v)\cup\{v\}\bigr)$, $N_G(X)=\{u\in V(G)\setminus X\mid u\text{ has a neighbor in }X\}$, and $M_G(X)=V(G)\setminus\bigl(X\cup N_G(X)\bigr)$. When there is no risk of confusion, we omit the subscript $G$ and simply write $N(v)$, $d(v)$, $M(v)$, $N(X)$, and $M(X)$. For $A,B\subseteq V(G)$, let $N_A(B)=N(B)\cap A$ and $M_A(B)=A\setminus\bigl(N_A(B)\cup B\bigr)$.

 	For a positive integer $k$, we use $[k]$ to denote the set $\{1,\ldots,k\}$. A {\em $k$-coloring} of a graph $G=(V,E)$ is a mapping $f:V\to [k]$ such that $f(u)\neq f(v)$ whenever $uv\in E$. We say that $G$ is {\em $k$-colorable} if it admits a $k$-coloring. The {\em chromatic number} of $G$, denoted by $\chi(G)$, is the smallest positive integer $k$ for which $G$ is $k$-colorable. A {\em clique} (resp. {\em stable set}) of $G$ is a set of pairwise adjacent (resp. nonadjacent) vertices of $G$. The {\em clique number} of $G$, denoted by $\omega(G)$, is the maximum cardinality of a clique in $G$.
 	
 	The concept of binding functions was introduced by Gy\'{a}rf\'{a}s \cite{G75} in 1975. Let ${\cal F}$ be a family of graphs. We say that ${\cal F}$ is $\chi$-{\em bounded} if there exists a function $f$ such that $\chi(H)\leqslant f(\omega(H))$ for every induced subgraph $H$ of every graph in ${\cal F}$. Such a function $f$ is called a {\em binding function} for ${\cal F}$.
 	
 	An induced cycle of length $k\geqslant 4$ is called a {\em hole}, and $k$ is called the {\em length} of the hole. A {\em $k$-hole} is a hole of length $k$. A hole is {\em odd} if its length is odd, and {\em even} otherwise. An {\em antihole} is the complement of a hole. A graph $G$ is said to be {\em perfect} if $\chi(H)=\omega(H)$ for every induced subgraph $H$ of $G$. The famous Strong Perfect Graph Theorem \cite{CRST06}, established by Chudnovsky {\em et al.} in 2006, states the following.
 	
 	\begin{theorem}\label{perfect}{\em \cite{CRST06}}
 		A graph $G$ is perfect if and only if $G$ is $(\text{odd hole},\text{odd antihole})$-free.
 	\end{theorem}
 	
 	\begin{figure}[htbp]
 		\begin{center}
 			\includegraphics[width=16cm]{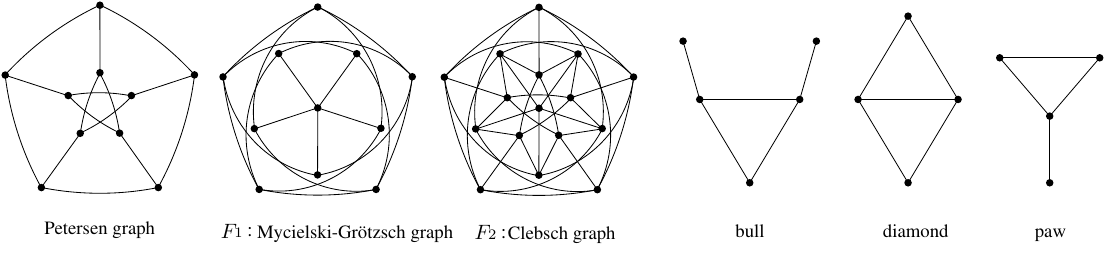}
 		\end{center}
 		\vskip -25pt
 		\caption{Illustration of the Petersen graph, $F_1$, $F_2$, the bull, the diamond, and the paw.}
 		\label{fig-1}
 	\end{figure}
 	
 	A graph $G$ is {\em $k$-vertex-critical} if $\chi(G)=k$ and every proper induced subgraph $H$ of $G$ satisfies $\chi(H)<k$. Let $F_1$ and $F_2$ denote the Mycielski-Gr\"{o}tzsch graph (Mycielski graph $G_4$) and the Clebsch graph, respectively. Notice that $F_1$ is an induced subgraph of $F_2$ and that $F_1$ is $4$-vertex-critical. A {\em bull} is a graph consisting of a triangle with two disjoint pendant edges; a {\em diamond} is a graph consisting of two triangles sharing exactly one edge; and a {\em paw} is a graph obtained from a triangle by adding a pendant edge (see Figure~\ref{fig-1}).

\begin{figure}[htbp]
 		\begin{center}
 			\includegraphics[width=14cm]{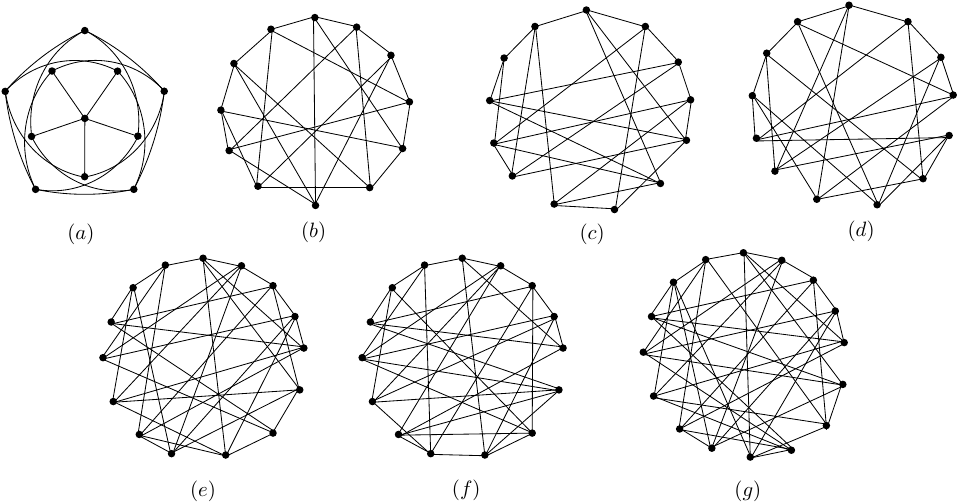}
 		\end{center}
 		\vskip -25pt
 		\caption{Illustration of (a)-(g).}
 		\label{fig-2}
 	\end{figure}

 	Let ${\cal G}$ be the family of seven explicitly defined graphs shown in Figure~\ref{fig-2}. Notice that graph (a) is isomorphic to $F_1$, and that every graph in ${\cal G}$ is $4$-vertex-critical and $(P_7,C_3)$-free \cite{G}. Using a computer-assisted approach, Goedgebeur and Schaudt \cite{G} determined all $4$-vertex-critical $(P_7,C_4)$-free graphs and proved that every $4$-vertex-critical $(P_7,C_3)$-free graph with at most 35 vertices belongs to ${\cal G}$. Moreover, they proposed the following conjecture.
 	
 	\begin{conjecture}\label{conj}
 		All $4$-vertex-critical $(P_7,C_3)$-free graphs belong to ${\cal G}$.
 	\end{conjecture}
 	
 	Zhou {\em et al.} \cite{Z} showed that Conjecture~\ref{conj} holds for $(P_7,C_3)$-free graphs containing $C_7$. In this paper, we prove that Conjecture~\ref{conj} also holds for $(P_2\cup P_4,C_3)$-free graphs. Our result is stated in the following theorem.
 	
 	\begin{theorem}\label{co-2}
 		The graph $F_1$ is the unique $4$-vertex-critical graph among the $(P_2\cup P_4,C_3)$-free graphs.
 	\end{theorem}
 	
 	Let $G$ be a graph. For a pair of nonadjacent vertices $u$ and $v$, we call $(u,v)$ a {\em comparable pair} if $N(u)\subseteq N(v)$. Note that if $(u,v)$ is a comparable pair in $G$, then $\chi(G)=\chi(G-u)$. We say that a graph $G$ is obtained from a graph $H$ by a {\em replication-vertex-addition} if $G$ is obtained by adding a vertex $u$ to $H$ such that there exists a vertex $v\in V(H)$ for which $(u,v)$ is a comparable pair in $G$. In fact, Theorem~\ref{co-2} follows directly from the structural theorem below, which is stronger than Theorem~\ref{co-2}.
 	
 	\begin{theorem}\label{structure}
 		Let $G$ be a connected $(P_2\cup P_4,C_3)$-free graph. Then one of the following holds.
 		\begin{enumerate}[label=(\roman*)]
 			\item $G$ has a comparable pair;
 			\item $\chi(G)\leqslant 3$;
 			\item $G$ contains $F_1$ as an induced subgraph and is an induced subgraph of $F_2$. Consequently, $\chi(G)=4$.
 		\end{enumerate}
 	\end{theorem}

 	In \cite{RST2002}, Randerath {\em et al.} proved that every $(P_6,C_3)$-free graph $G$ satisfies $\chi(G)\leqslant 4$, and that every such graph with $\chi(G)=4$ contains $F_1$. They also gave a polynomial-time algorithm for deciding whether a $(P_6,C_3)$-free graph is $3$-colorable. In \cite{Z}, Zhou {\em et al.} proved that every $(P_5\cup P_1,C_3)$-free graph $G$ satisfies $\chi(G)\leqslant 3$. In \cite{P2013}, Pyatkin proved that every $(2P_3,C_3)$-free graph $G$ satisfies $\chi(G)\leqslant 4$. In \cite{BC2018}, Bharathi and Choudum proved that every $(P_2\cup P_4)$-free graph $G$ satisfies $\chi(G)\leqslant\binom{\omega(G)+2}{3}$, which implies that every $(P_2\cup P_4,C_3)$-free graph is $4$-colorable.
 	
 	By Theorem~\ref{structure}, we can determine the chromatic number of $(P_2\cup P_4,C_3)$-free graphs. In particular, every $(P_2\cup P_4,C_3)$-free graph is $4$-colorable, which can also be shown by a simple induction on $|V(G)|$. Moreover, the proof of Theorem~\ref{structure} yields a polynomial-time algorithm for deciding whether a $(P_2\cup P_4,C_3)$-free graph is $3$-colorable. The following corollary follows immediately from Theorems~\ref{perfect} and~\ref{structure}.

 	\begin{corollary}\label{co-1}
 		Let $G$ be a connected $(P_2\cup P_4,C_3)$-free graph. Then the following hold.
 		\begin{enumerate}[label=(\roman*)]
 			\item $\chi(G)=4$ if and only if $G$ can be obtained from a graph $H$, which contains $F_1$ as an induced subgraph and is an induced subgraph of $F_2$, by a sequence of replication-vertex-additions;
 			\item $\chi(G)=3$ if and only if $G$ contains a $5$-hole or a $7$-hole and is $F_1$-free.
 		\end{enumerate}
 	\end{corollary}
 	
 	In \cite{O1988}, Olariu showed that every connected paw-free graph is either triangle-free or complete multipartite. Hence, Theorem~\ref{structure} immediately yields the following corollary.
 	
 	\begin{corollary}\label{co-3}
 		Let $G$ be a connected $(P_2\cup P_4,paw)$-free graph. Then one of the following holds.
 		\begin{enumerate}[label=(\roman*)]
 			\item $G$ has a comparable pair;
 			\item $G$ is a complete graph;
 			\item $G$ contains $F_1$ as an induced subgraph and is an induced subgraph of $F_2$. Consequently, $\chi(G)=4$;
 			\item $\chi(G)\leqslant 3$, and there is a polynomial-time algorithm for finding a $3$-coloring of $G$.
 		\end{enumerate}
 	\end{corollary}

 	A graph $G$ is said to be {\em perfectly divisible} if, for every induced subgraph $H$ of $G$, the vertex set $V(H)$ can be partitioned into two sets $A$ and $B$ such that $H[A]$ is perfect and $\omega(H[B])<\omega(H)$. This concept was introduced by Ho\'ang \cite{H22}. Chudnovsky and Sivaraman \cite{CS19} proved that every $(P_5,\text{bull})$-free graph and every $(\text{odd hole},\text{bull})$-free graph is perfectly divisible. Chen and Xu \cite{CX2025} proved that every $(P_7,C_5,\text{bull})$-free graph is perfectly divisible.
 	
 	Let $S\subseteq V(G)$ with $1<|S|<|V(G)|$. We say that $S$ is a {\em homogeneous set} of $G$ if every vertex in $V(G)\setminus S$ is either complete or anticomplete to $S$.
 	
 	Notice that $F_1$ is a $(P_2\cup P_4,\text{bull})$-free graph that is not perfectly divisible, and $\omega(F_1)=2$. Thus, not every $(P_2\cup P_4,\text{bull})$-free graph is perfectly divisible. Recently, Deng and Chang \cite{DC2025} proved that every $(P_2\cup P_3,\text{bull})$-free graph $G$ with $\omega(G)\geqslant 3$ and no homogeneous set admits a partition $(X,Y)$ such that $G[X]$ is perfect and $\omega(G[Y])<\omega(G)$. Subsequently, Chen and Wang \cite{CW2025} extended this result to the larger class of $(P_2\cup P_4,\text{bull})$-free graphs. However, the result does not hold for all $(P_2\cup P_4,\text{bull})$-free graphs. Indeed, let $G$ be obtained from $F_1$ by adding a copy of $K_n$ such that every vertex of the $K_n$ is complete to $V(F_1)$. Then $G$ is $(P_2\cup P_4,\text{bull})$-free, $\omega(G)=n+2$, and $G$ is not perfectly divisible.
 	
 	In this paper, we prove the following theorem, which generalizes the main results of Deng and Chang \cite{DC2025} and Chen and Wang \cite{CW2025}. 
 	
 	\begin{theorem}\label{perfectdivisible}
 		Let $G$ be a $(P_2\cup P_4,\text{bull})$-free graph. Then $G$ is perfectly divisible if and only if $G$ is $F_1$-free.
 	\end{theorem}
 	
 	By a simple induction on $\omega(G)$, every perfectly divisible graph $G$ satisfies $\chi(G)\leqslant \binom{\omega(G)+1}{2}$ \cite{CS19}. Therefore, Theorem~\ref{perfectdivisible} immediately yields the following corollary.
 	
 	\begin{corollary}\label{chromatic}
 		Let $G$ be a $(P_2\cup P_4,\text{bull},F_1)$-free graph. Then $\chi(G)\leqslant \binom{\omega(G)+1}{2}$.
 	\end{corollary}
 	
 	Notice that the class of $3K_1$-free graphs does not admit a linear binding function \cite{BGS2002,SR2019}. The same is true for the class of $(P_2\cup P_4,\text{bull},F_1)$-free graphs.
 	
 	Denote the {\em minimum (resp. maximum) degree} of $G$ by $\delta(G)$ (resp. $\Delta(G)$). For $u,v\in V(G)$, we write $u\sim v$ if $uv\in E(G)$, and $u\not\sim v$ if $uv\notin E(G)$.
 	
 	The {\em Cartesian product} of two graphs $G$ and $H$, denoted by $G\square H$, is the graph with vertex set $\{(x,y)\mid x\in V(G),y\in V(H)\}$, where two vertices $(x_1,y_1)$ and $(x_2,y_2)$ are adjacent in $G\square H$ if and only if either $x_1=x_2$ and $y_1\sim y_2$ in $H$, or $y_1=y_2$ and $x_1\sim x_2$ in $G$. In \cite{CX2024}, Chen and Xu proved that if $G$ is a connected $(\text{bull},\text{diamond})$-free graph with $\omega(G)\geqslant3$, then either $\delta(G)\leqslant\omega(G)-1$ or $G$ is isomorphic to $K_2\square K_{\omega(G)}$. Theorem~\ref{structure} yields the following corollary.
 	
 	\begin{corollary}\label{diamond}
 		Let $G$ be a connected $(P_2\cup P_4,\text{bull},\text{diamond})$-free graph. Then one of the following holds.
 		\begin{enumerate}[label=(\roman*)]
 			\item $G$ has a comparable pair;
 			\item $\delta(G)\leqslant\omega(G)-1$;
 			\item $G$ is isomorphic to $K_2\square K_{\omega(G)}$;
 			\item $G$ contains $F_1$ as an induced subgraph and is an induced subgraph of $F_2$. Consequently, $\chi(G)=4$;
 			\item $\chi(G)\leqslant3$, and there is a polynomial-time algorithm for finding a $3$-coloring of $G$.
 		\end{enumerate}
 	\end{corollary}
 	
 	According to Corollary \ref{diamond}, by a simple induction on $|V(G)|$, we can derive the following corollary. 
 	
 	\begin{corollary}\label{co-4}
 		Let $G$ be a $(P_2\cup P_4,\text{bull},\text{diamond})$-free graph. Then $\chi(G)\leqslant \max\{4,\omega(G)\}$. 
 	\end{corollary}
 	
 	Angeliya {\em et al.} \cite{AKH2025} proved that every $(P_2\cup P_4,\text{diamond})$-free graph $G$ satisfies $\chi(G)\leqslant\max\{6,\omega(G)\}$. The bound in Corollary~\ref{co-4} is optimal and extends their result to bull-free graphs.
 	
 	The remainder of the paper is organized as follows. In Section~2, we prove Theorem~\ref{structure} by analyzing the structure of $(P_2\cup P_4,C_3)$-free graphs. In Section~3, we use this structural theorem to obtain an $O(n^4)$-time algorithm for deciding $3$-colorability of $(P_2\cup P_4,C_3)$-free graphs. In Section~4, we prove Theorem~\ref{perfectdivisible}.

 	\section{Proof of Theorem~\ref{structure}}\label{1.2}
 	
 	In this section, we prove Theorem~\ref{structure}. Brooks \cite{B41} showed that for a graph $G$, the chromatic number is at most $\Delta(G)$, unless $G$ is a complete graph or an odd cycle. In either of these two exceptional cases, $\chi(G)=\Delta(G)+1$.
 	
 	\begin{theorem}\label{Brook}{\em \cite{B41}}
 		Let $G$ be a graph with $\Delta(G)\geqslant 3$. Then $\chi(G)\leqslant\max\{\Delta(G),\omega(G)\}$.
 	\end{theorem}
 	
 	A {\em dominating set} in a graph $G$ is a subset $S$ of $V(G)$ such that every vertex in $V(G)\setminus S$ has a neighbor in $S$. We first present the following two lemmas.
 	
 	\begin{lemma}\label{c5}
 		Let $G$ be a $(P_2\cup P_4,\text{bull})$-free graph, let $v\in V(G)$, and let $C=v_1v_2v_3v_4v_5v_1$ be a $5$-hole in $G[M(v)]$. Then, for every vertex $x\in N(v)$, either $N(x)\cap V(C)=\{v_i,v_{i+2}\}$ for some $i\in[5]$, or $N(x)\cap V(C)=V(C)$. (The subscript is taken modulo $5$.)
 	\end{lemma}
 	
 	\pf
 	Suppose, to the contrary, that there exists $x\in N(v)$ for which neither of the two conclusions holds. To avoid an induced $P_2\cup P_4$ on $\{v,x\}\cup V(C)$, we have $N_C(x)\neq\emptyset$. Without loss of generality, we may assume that $x\sim v_1$. Suppose that $N_C(x)\neq V(C)$. Then $x$ has a neighbor in $\{v_2,v_3,v_4,v_5\}$, as otherwise $\{x,v,v_2,v_3,v_4,v_5\}$ induces a $P_2\cup P_4$.
 	
 	If $x\sim v_2$, then $x\sim v_3$ to avoid an induced bull on $\{x,v,v_1,v_2,v_3\}$. Similarly, $x\sim v_5$. In this case, we have $x\sim v_4$, as otherwise $\{x,v,v_2,v_3,v_4\}$ induces a bull. Hence, $x$ is complete to $V(C)$, a contradiction. Therefore, $x\not\sim v_2$, and similarly, $x\not\sim v_5$.
 	
 	Now, $x$ has a neighbor in $\{v_3,v_4\}$. If $x$ is complete to $\{v_3,v_4\}$, then $\{x,v,v_2,v_3,v_4\}$ induces a bull, a contradiction. Thus, $x$ has exactly one neighbor in $\{v_3,v_4\}$. Hence, $N_C(x)=\{v_1,v_3\}$ or $\{v_4,v_1\}$, contradicting our assumption. This proves Lemma~\ref{c5}. \qed
 	
 	\begin{lemma}\label{c7}
 		Let $G$ be a connected $(P_2\cup P_4,\text{bull})$-free graph, and let $C=v_1v_2v_3v_4v_5v_6v_7v_1$ be a $7$-hole in $G$. If there does not exist a vertex that is complete to $V(C)$, then $V(C)$ is a dominating set of $G$. (The subscript is taken modulo $7$.)
 	\end{lemma}
 	
 	\pf
 	Suppose that there does not exist a vertex that is complete to $V(C)$. For $i\in[7]$, let
 	\begin{eqnarray*}
 		X_i&=&\{u\in N(V(C))\mid N_C(u)=\{v_i,v_{i+2}\}\};\\
 		Y_i&=&\{u\in N(V(C))\mid N_C(u)=\{v_i,v_{i+2},v_{i+4}\}\}.
 	\end{eqnarray*}
 	
 	Let $X=\bigcup_{i=1}^7X_i$ and $Y=\bigcup_{i=1}^7Y_i$. We next prove the following claim.
 	
 	\begin{claim}\label{c-1}
 		$N(V(C))=X\cup Y$.
 	\end{claim}
 	
 	\pf
 	It suffices to prove that $N(V(C))\subseteq X\cup Y$. Let $x\in N(V(C))$. Without loss of generality, suppose that $x\sim v_1$. To avoid an induced $P_2\cup P_4$ on $\{x,v_1,v_3,v_4,v_5,v_6\}$, $x$ has a neighbor in $\{v_3,v_4,v_5,v_6\}$.
 	
 	Suppose that $x$ is anticomplete to $\{v_3,v_6\}$. Then $x$ has a neighbor in $\{v_4,v_5\}$. Moreover, $x$ has exactly one neighbor in $\{v_4,v_5\}$, as otherwise $\{x,v_3,v_4,v_5,v_6\}$ induces a bull. Without loss of generality, suppose that $x\sim v_4$ and $x\not\sim v_5$. If $x\sim v_7$, then $\{x,v_1,v_4,v_6,v_7\}$ induces a bull. Thus, $x\not\sim v_7$, and then $x\not\sim v_2$ to avoid a bull on $\{x,v_1,v_2,v_3,v_7\}$. But now, $\{x,v_2,v_3,v_4,v_6,v_7\}$ induces a $P_2\cup P_4$. Therefore, $x$ has a neighbor in $\{v_3,v_6\}$. Without loss of generality, we may assume that $x\sim v_3$.
 	
 	Suppose that $x\sim v_2$. If $x\sim v_4$, then $x\sim v_5$ to avoid a bull on $\{x,v_1,v_3,v_4,v_5\}$. Consequently, $x\sim v_6$ to avoid a bull on $\{x,v_1,v_4,v_5,v_6\}$. Since $x$ is not complete to $V(C)$ by our assumption, it follows that $\{x,v_2,v_5,v_6,v_7\}$ induces a bull. Therefore, $x\not\sim v_2$.
 	
 	Suppose that $x\sim v_4$. Then $x\sim v_5$, as otherwise $\{x,v_2,v_3,v_4,v_5\}$ induces a bull. But now, $\{x,v_1,v_4,v_5,v_6\}$ induces a bull if $x\not\sim v_6$, and $\{x,v_2,v_3,v_4,v_6\}$ induces a bull if $x\sim v_6$. Both are contradictions. Hence, $x\not\sim v_4$. Similarly, $x\not\sim v_7$.
 	
 	If $x$ is complete to $\{v_5,v_6\}$, then $\{x,v_4,v_5,v_6,v_7\}$ induces a bull, a contradiction. Therefore, $N_C(x)\in\{\{v_1,v_3\},\{v_1,v_3,v_5\},\{v_1,v_3,v_6\}\}$. Thus, $x\in X_1\cup Y_1\cup Y_6$. This proves Claim~\ref{c-1}. \qed
 	
 	\medskip
 	
 	Recall that $M(V(C))=V(G)\setminus(N(V(C))\cup V(C))$. To prove that $V(C)$ is a dominating set of $G$, it suffices to show that $M(V(C))=\emptyset$. Suppose, to the contrary, that $M(V(C))\neq\emptyset$. Since $G$ is connected, there exist two vertices $u,v\in V(G)$ such that $v\in M(V(C))$, $u\in N(V(C))$, and $u\sim v$. By Claim~\ref{c-1}, we may assume by symmetry that $u\in X_1\cup Y_1$. But then, $\{u,v,v_2,v_3,v_6,v_7\}$ induces a $P_2\cup P_4$, a contradiction. This proves Lemma~\ref{c7}. \qed
 	
 	\medskip
 	
 	Now, we proceed to prove Theorem~\ref{structure}.
 	
 	\noindent\textbf{{\em Proof of Theorem~\ref{structure}:}}
 	Let $G$ be a connected $(P_2\cup P_4,C_3)$-free graph. Suppose that $G$ has no comparable pair and $\chi(G)\geqslant4$. Since $G$ is neither an odd cycle nor a complete graph, by Theorem~\ref{Brook}, we have $\Delta(G)\geqslant4$. Let $v\in V(G)$ with $d(v)=\Delta(G)$ and let $G'=G[M(v)]$. Notice that $N(v)$ is a stable set since $G$ is triangle-free. Thus, $G'$ is not bipartite, as otherwise $\chi(G)\leqslant3$, a contradiction.
 	
 	Since $G'$ is a $(P_2\cup P_4,C_3)$-free graph, every odd hole in $G'$ has length $5$ or $7$. Indeed, an odd hole of length at least $9$ contains an induced $P_2\cup P_4$. Moreover, $G'$ contains no odd antihole of length at least $7$, since the complement of such an antihole contains a triangle. Thus, by Theorem~\ref{perfect}, $G'$ contains a $5$-hole or a $7$-hole.
 	
 	We claim that $G'$ contains no $7$-hole. Suppose, to the contrary, that $C'$ is a $7$-hole in $G'$. Since $G$ is triangle-free, no vertex of $G$ is complete to $V(C')$. By Lemma~\ref{c7}, $V(C')$ is a dominating set of $G$. However, $C'\subseteq M(v)$, and hence $v$ is anticomplete to $V(C')$, a contradiction. Therefore, $G'$ contains a $5$-hole $C=v_1v_2v_3v_4v_5v_1$. From now on, the subscripts are taken modulo $5$ in the proof of Theorem~\ref{structure}. We begin with the following claim.
 	
 	\begin{claim}\label{claim1}
 		Let $u\in N(v)$. Then $N_C(u)=\{v_i,v_{i+2}\}$ for some $i\in[5]$.
 	\end{claim}
 	
 	\pf
 	Since $G$ is triangle-free, we have $N_C(u)\neq V(C)$. By Lemma~\ref{c5}, $N_C(u)=\{v_i,v_{i+2}\}$ for some $i\in[5]$. This proves Claim~\ref{claim1}. \qed
 	
 	\begin{claim}\label{claim3}
 		$G'$ is connected.
 	\end{claim}
 	
 	\pf
 	Suppose, to the contrary, that there exists a component $T$ of $G'$ different from the one containing $C$. Since $G$ is connected, there exist vertices $u\in V(T)$ and $w\in N(v)$ such that $u\sim w$. Without loss of generality, suppose that $N_C(w)=\{v_1,v_3\}$ by Claim~\ref{claim1}. Notice that $u\not\sim v$ and $w\in N(v)\cap N(u)$. Since $G$ has no comparable pair, there exists a vertex $u'\in N(u)\setminus N(v)$. Clearly, $u'\in V(T)$, and thus $u'$ is anticomplete to $V(C)$. But then, $\{u,u',v_1,v_2,v_3,v_4\}$ induces a $P_2\cup P_4$, a contradiction. This proves Claim~\ref{claim3}. \qed
 	
 	\begin{claim}\label{claim2}
 		For each $i\in[5]$, there is at most one vertex in $N(v)$ that is complete to $\{v_i,v_{i+2}\}$, and hence $4\leqslant\Delta(G)\leqslant5$.
 	\end{claim}
 	
 	\pf
 	Without loss of generality, set $i=1$. Suppose that there exist two vertices $w_1,w_2\in N(v)$ such that $\{w_1,w_2\}$ is complete to $\{v_1,v_3\}$. By Claim~\ref{claim1}, $N_C(w_1)=N_C(w_2)=\{v_1,v_3\}$. Notice that $\{v,v_1,v_3\}\subseteq N(w_1)\cap N(w_2)$ and $w_1\not\sim w_2$. Since $G$ has no comparable pair and is triangle-free, there exists a vertex $x\in V(G)\setminus(V(C)\cup N(v)\cup\{v\})$ such that $x\sim w_2$ and $x\not\sim w_1$. Moreover, $x$ is anticomplete to $\{v_1,v_3\}$ to avoid triangles. To avoid an induced $P_2\cup P_4$ on $\{v_4,v_5,w_1,v,w_2,x\}$, we have that either $x\sim v_4$ or $x\sim v_5$.
 	
 	Suppose that $x\sim v_4$. Then $x\sim v_2$, as otherwise $\{x,v_4,v,w_1,v_1,v_2\}$ induces a $P_2\cup P_4$. Thus, $N_C(x)=\{v_2,v_4\}$ since $G$ is triangle-free. But then, $\{v,w_1,v_2,x,v_4,v_5\}$ induces a $P_2\cup P_4$, a contradiction. Therefore, $x\not\sim v_4$, and hence $x\sim v_5$.
 	
 	To avoid an induced $P_2\cup P_4$ on $\{x,v_5,v,w_1,v_3,v_2\}$, we have that $x\sim v_2$. But now, $\{w_1,v,v_2,x,v_5,v_4\}$ induces a $P_2\cup P_4$, a contradiction. This proves that for each $i\in[5]$, there is at most one vertex in $N(v)$ that is complete to $\{v_i,v_{i+2}\}$. Thus, $\Delta(G)\leqslant5$ by Claim~\ref{claim1}. Since $\Delta(G)\geqslant4$, we conclude that $4\leqslant\Delta(G)\leqslant5$. This proves Claim~\ref{claim2}. \qed
 	
 	\begin{claim}\label{c-11}
 		$\Delta(G)=5$.
 	\end{claim}
 	
 	\pf
 	Suppose, to the contrary, that $\Delta(G)=4$. By Claims~\ref{claim1} and~\ref{claim2}, we may assume by symmetry that $N(v)=\{w_1,w_2,w_3,w_4\}$ and $N_C(w_i)=\{v_i,v_{i+2}\}$ for $i\in[4]$. Since $\chi(G)\geqslant4$, it follows that $V(G')\setminus V(C)\neq\emptyset$. Let $Y=N_{G'}(V(C))$. By Claim~\ref{claim3}, $G'$ is connected, and hence $Y\neq\emptyset$. Moreover, we have $d(v_1)=d(v_3)=d(v_4)=4=\Delta(G)$, and thus
 	\begin{equation}\label{e-11}
 		\mbox{every vertex in $Y$ is adjacent to $v_2$ or $v_5$.}
 	\end{equation}
 	
 	We next prove that
 	\begin{equation}\label{e-1}
 		\mbox{for every vertex $y\in Y$, $y$ is not complete to $\{v_2,v_5\}$.}
 	\end{equation}
 	
 	Suppose, to the contrary, that $y$ is complete to $\{v_2,v_5\}$. To avoid an induced $P_2\cup P_4$ on $\{v,w_1,v_2,y,v_5,v_4\}$ or $\{v,w_4,v_3,v_2,y,v_5\}$, $y$ must be complete to $\{w_1,w_4\}$. Then $d(y)=4=\Delta(G)$, and thus $N(y)=N(v_1)$. It follows that $(y,v_1)$ is a comparable pair of $G$, a contradiction. This proves (\ref{e-1}).
 	
 	By (\ref{e-11}) and (\ref{e-1}), there exists a partition $(Y_2,Y_5)$ of $Y$ such that each vertex in $Y_2$ belongs to $N(v_2)\setminus N(v_5)$ and each vertex in $Y_5$ belongs to $N(v_5)\setminus N(v_2)$. Since $\Delta(G)=4$, for each vertex $y_2\in Y_2$, we have $N_C(y_2)=\{v_2\}$ and $|Y_2|\leqslant1$, and for each vertex $y_5\in Y_5$, we have $N_C(y_5)=\{v_5\}$ and $|Y_5|\leqslant1$. Moreover, by Lemma~\ref{c5}, every vertex in $Y$ has no neighbor in $M(V(C))$. Since $G'$ is connected, it follows that $M(V(C))=\emptyset$, and hence $V(G)=V(C)\cup N(v)\cup Y_2\cup Y_5\cup\{v\}$.
 	
 	Now, we may construct a $3$-coloring $\phi$ of $G$ by setting $\phi(\{v,v_1,v_3\}\cup Y_2)=1$, $\phi(\{v_2,v_4,w_3\}\cup Y_5)=2$, and $\phi(\{v_5,w_1,w_2,w_4\})=3$, a contradiction. This proves Claim~\ref{c-11}. \qed
 	
 	\medskip
 	
 	By Claim~\ref{c-11}, we have $\Delta(G)=5$. Without loss of generality, we may suppose that $N(v)=\{w_1,w_2,w_3,w_4,w_5\}$ and $N_C(w_i)=\{v_i,v_{i+2}\}$ for each $i\in[5]$ by Claims~\ref{claim1} and~\ref{claim2}. Then $G$ contains $F_1$, since $G[N(v)\cup\{v\}\cup V(C)]$ is isomorphic to $F_1$.
 	
 	\begin{claim}\label{5y}
 		For each vertex $y\in V(G')\setminus V(C)$, if $N(y)\cap V(C)\neq\emptyset$, then $N(y)\cap V(C)=\{v_i\}$ for some $i\in[5]$.
 	\end{claim}
 	
 	\pf
 	Suppose, to the contrary, that there exists a vertex $y\in V(G')\setminus V(C)$ such that $N(y)\cap V(C)\neq\emptyset$ and $N(y)\cap V(C)\neq\{v_i\}$ for each $i\in[5]$. Since $G$ is triangle-free, we have $N(y)\cap V(C)=\{v_i,v_{i+2}\}$. Without loss of generality, set $i=1$. Then $y\sim w_5$, as otherwise $\{v,w_5,v_1,y,v_3,v_4\}$ induces a $P_2\cup P_4$. Similarly, to avoid an induced $P_2\cup P_4$ on $\{v,w_2,v_3,y,v_1,v_5\}$, we have $y\sim w_2$. Since $G$ is triangle-free, $y$ is anticomplete to $\{w_1,w_3,w_4\}$.
 	
 	Notice that $\{v_1,v_3,w_2,w_5\}\subseteq N(v_2)\cap N(y)$ and $v_2\not\sim y$. Since $G$ has no comparable pair, it follows that $N(y)\not\subseteq N(v_2)$, and thus there exists a vertex $y'$ such that $y'\sim y$ and $y'\not\sim v_2$. Clearly, $y'$ is anticomplete to $\{v_1,v_3,w_2,w_5\}$ since $G$ is triangle-free. To avoid an induced $P_2\cup P_4$ on $\{v,w_4,v_2,v_3,y,y'\}$, we have $y'\sim w_4$, and hence $y'\not\sim v_4$ since $G$ is triangle-free. Therefore,
 	\begin{equation}\label{e1}
 		\mbox{$y'$ is anticomplete to $\{v_1,v_2,v_3,v_4,w_2,w_5\}$ and $y'\sim w_4$.}
 	\end{equation}
 	
 	To avoid an induced $P_2\cup P_4$ on $\{w_3,v_5,w_2,y,y',w_4\}$, $y'$ is adjacent to $w_3$ or $v_5$. Next, we prove that
 	\begin{equation}\label{w3}
 		y'\not\sim w_3.
 	\end{equation}
 	
 	Suppose, to the contrary, that $y'\sim w_3$. Then $y'\not\sim v_5$, as otherwise $y'w_3v_5y'$ is a triangle. Together with (\ref{e1}), this implies that $y'$ is anticomplete to $V(C)$. To avoid an induced $P_2\cup P_4$ on $\{y',w_4,w_1,v_3,v_2,w_5\}$, we have $y'\sim w_1$. But then, $\{y',w_1,v_2,w_5,v_5,v_4\}$ induces a $P_2\cup P_4$ by (\ref{e1}), a contradiction. This proves (\ref{w3}). \medskip
 	
 	By (\ref{w3}), we have $y'\sim v_5$, and $N_C(y')=\{v_5\}$ by (\ref{e1}). But then, $\{y',v_5,v,w_2,v_2,v_3\}$ induces a $P_2\cup P_4$, a contradiction. This proves Claim~\ref{5y}. \qed
 	
 	By Claim~\ref{claim3}, $G'$ is connected. Therefore, by Claim~\ref{5y} and Lemma~\ref{c5}, we deduce that $M_{G'}(V(C))=\emptyset$, and for every vertex $y\in V(G')\setminus V(C)$, there exists an $i\in[5]$ such that
 	\begin{equation}\label{yunique}
 		\mbox{$N_C(y)=\{v_i\}$.}
 	\end{equation}
 	
 	Furthermore, the condition $\Delta(G)=5$ implies that for each $i\in[5]$,
 	\begin{equation}\label{vunique}
 		v_i~\text{has at most one neighbor in}~V(G')\setminus V(C).
 	\end{equation}
 	
 	Let $Y_i=N_{G'}(v_i)$ for $i\in[5]$. By (\ref{yunique}) and (\ref{vunique}), we have $\bigcup_{i=1}^5Y_i=V(G')\setminus V(C)$, $|Y_i|\leqslant1$, and for any vertex $y_i\in Y_i$, $N_C(y_i)=\{v_i\}$. Moreover,
 	\begin{equation}\label{G}
 		V(G)=N(v)\cup\{v\}\cup V(C)\cup\left(\bigcup_{i=1}^5Y_i\right).
 	\end{equation}
 	Hence, $|V(G)|\leqslant16$.
 	
 	For each $i\in[5]$, since $|Y_i|\leqslant1$, we may assume that $Y_i=\{y_i\}$ whenever $Y_i\neq\emptyset$ in the remainder of the proof. Since $G$ is triangle-free, we have
 	\begin{equation}\label{y14}
 		y_i~\text{is anticomplete to}~\{w_i,w_{i+3}\}.
 	\end{equation}
 	
 	\begin{claim}\label{y1w2w3}
 		$N(y_i)\cap N(v)=\{w_{i+1},w_{i+2}\}$.
 	\end{claim}
 	
 	\pf
 	Without loss of generality, set $i=1$. To avoid an induced $P_2\cup P_4$ on $\{v_1,y_1,w_2,v_4,v_3,w_3\}$, we have $y_1\sim w_2$ or $y_1\sim w_3$. If $y_1\sim w_2$ and $y_1\not\sim w_3$, then $\{w_3,v_3,w_2,y_1,v_1,w_4\}$ induces a $P_2\cup P_4$, a contradiction. Conversely, if $y_1\sim w_3$ and $y_1\not\sim w_2$, then $\{w_2,v_4,w_3,y_1,v_1,w_1\}$ induces a $P_2\cup P_4$, also a contradiction. Therefore, $y_1$ is complete to $\{w_2,w_3\}$. Moreover, $y_1\not\sim w_5$, as otherwise $\{y_1,w_5,w_1,v_3,v_4,w_4\}$ induces a $P_2\cup P_4$ by (\ref{y14}). Hence, $N(y_1)\cap N(v)=\{w_2,w_3\}$. This proves Claim~\ref{y1w2w3}. \qed
 	
 	\medskip
 	
 	\begin{claim}\label{petersen}
 		$Y_i$ is anticomplete to $Y_{i+1}\cup Y_{i-1}$ and complete to $Y_{i+2}\cup Y_{i-2}$.
 	\end{claim}
 	
 	\pf
 	Without loss of generality, set $i=1$. Suppose, to the contrary, that $y_1\sim y_2$. By Claim~\ref{y1w2w3}, $w_3$ is complete to $\{y_1,y_2\}$, and then $y_1y_2w_3y_1$ is a triangle, a contradiction. Therefore, $Y_1$ is anticomplete to $Y_2\cup Y_5$ by symmetry.
 	
 	If $y_1\not\sim y_3$, then $\{y_1,w_2,v_3,y_3,w_5,v_5\}$ induces a $P_2\cup P_4$ by Claim~\ref{y1w2w3}, a contradiction. Thus, $Y_1$ is complete to $Y_3\cup Y_4$ by symmetry. This proves Claim~\ref{petersen}. \qed
 	
 	\medskip
 	
 	By Claims~\ref{y1w2w3} and~\ref{petersen}, we have that $G[(\bigcup_{i=1}^5Y_i)\cup V(C)]$ is an induced subgraph of the Petersen graph (see Figure~\ref{fig-1}). Therefore, by (\ref{G}), $G$ is an induced subgraph of $F_2$. Since $G[N(v)\cup\{v\}\cup V(C)]$ is isomorphic to $F_1$, $G$ contains $F_1$ as an induced subgraph. This completes the proof of Theorem~\ref{structure}. \qed
 	
 	\section{A polynomial-time algorithm for $3$-colorability}\label{algorithm}
 	
 	In this section, we show that Theorem~\ref{structure} yields a polynomial-time algorithm for deciding whether a $(P_2\cup P_4,C_3)$-free graph is $3$-colorable.
 	
 	Recall that if $(u,v)$ is a comparable pair in a graph $G$, then $u\not\sim v$ and $N(u)\subseteq N(v)$. In this case, $\chi(G)=\chi(G-u)$. Therefore, deleting $u$ does not change whether $G$ is $3$-colorable.
 	
 	Let $G$ be a $(P_2\cup P_4,C_3)$-free graph with $n$ vertices. Since $G$ is $3$-colorable if and only if every connected component of $G$ is $3$-colorable, it suffices to consider each connected component separately.
 	
 	Starting with a connected component of $G$, whenever the current graph contains a comparable pair $(u,v)$, we delete $u$. We repeat this process until no comparable pair remains, and let $H$ be the resulting graph. Since each deletion preserves the chromatic number, the original component is $3$-colorable if and only if $H$ is $3$-colorable.
 	
 	A comparable pair can be found in $O(n^3)$ time. Indeed, we examine every pair of nonadjacent vertices $u$ and $v$ and check whether $N(u)\subseteq N(v)$ or $N(v)\subseteq N(u)$. There are $O(n^2)$ pairs of vertices, and each neighborhood inclusion can be checked in $O(n)$ time. Since at most $n$ vertices can be deleted, the entire reduction process takes $O(n^4)$ time.
 	
 	Now $H$ has no comparable pair. By Theorem~\ref{structure}, one of the following holds:
 	\begin{enumerate}[label=(\roman*)]
 		\item $\chi(H)\leqslant3$;
 		\item $H$ contains $F_1$ as an induced subgraph and is an induced subgraph of $F_2$.
 	\end{enumerate}
 	
 	Moreover, the proof of Theorem~\ref{structure} shows that in the second case, $|V(H)|\leqslant16$. Therefore, if $|V(H)|>16$, then necessarily $\chi(H)\leqslant3$. If $|V(H)|\leqslant16$, we can determine whether $H$ is $3$-colorable by checking all possible $3$-colorings of $H$. Since there are at most $3^{16}$ such colorings, this step takes constant time.
 	
 	Consequently, the entire algorithm runs in $O(n^4)$ time. Thus, there is an $O(n^4)$-time algorithm for deciding whether a $(P_2\cup P_4,C_3)$-free graph with $n$ vertices is $3$-colorable.

 	\section{Proof of Theorem~\ref{perfectdivisible}}\label{1.3}
 	
 	In this section, we prove Theorem~\ref{perfectdivisible}. Let $G$ be a graph and let $w$ be a weight function on $V(G)$. We use {\em $\omega_w(G)$} to denote the maximum weight of a clique in $G$. A graph $G$ is {\em perfectly weight divisible} if, for every positive integer weight function $w$ on $V(G)$ and every induced subgraph $H$ of $G$, there is a partition $(P,W)$ of $V(H)$ such that $H[P]$ is perfect and $\omega_w(H[W])<\omega_w(H)$. Clearly, every perfectly weight divisible graph is perfectly divisible.
 	
 	We will use the following lemmas.
 	
 	\begin{lemma}\label{homogeneous}{\em \cite{CS19}}
 		Every minimal nonperfectly weight divisible graph has no homogeneous set.
 	\end{lemma}
 	
 	\begin{lemma}\label{3-coloring}
 		Let $G$ be a graph with $\chi(G)\leqslant3$. Then $G$ is perfectly weight divisible.
 	\end{lemma}
 	
 	\pf
 	Let $w$ be a positive integer weight function on $V(G)$, and let $G'$ be an induced subgraph of $G$. Since $\chi(G')\leqslant3$, there is a partition $(A,B)$ of $V(G')$ such that $G'[A]$ is bipartite and $B$ is a stable set of $G'$. Let $X\subseteq B$ be the set of isolated vertices of $G'$ contained in $B$. Then $(A\cup X,B\setminus X)$ is a partition of $V(G')$, and $G'[A\cup X]$ is perfect.
 	
 	If $B\setminus X=\emptyset$, then $\omega_w(G'[B\setminus X])=0<\omega_w(G')$. Otherwise, since $B\setminus X$ is a stable set, let $b\in B\setminus X$ be a vertex of maximum weight in $B\setminus X$. Since $b$ is not isolated in $G'$, it has a neighbor $a\in A$. As $w(a)>0$, we have $\omega_w(G')\geqslant w(a)+w(b)>w(b)=\omega_w(G'[B\setminus X])$. Hence, $G$ is perfectly weight divisible. This proves the lemma. \qed
 	
 	\begin{lemma}\label{bullfree}{\em \cite{CS2008}}
 		If $G$ is a bull-free graph, then either $G$ has a homogeneous set or, for every $v\in V(G)$, either $G[N(v)]$ is perfect or $G[M(v)]$ is perfect.
 	\end{lemma}
 	
 	\noindent\textbf{{\em Proof of Theorem~\ref{perfectdivisible}:}}
 	Let $G$ be a $(P_2\cup P_4,\text{bull})$-free graph. First, suppose that $G$ is perfectly divisible. Since $F_1$ is not perfectly divisible, it follows that $G$ is $F_1$-free.
 	
 	Now, suppose that $G$ is $F_1$-free. To prove the converse, it suffices to show that every $(P_2\cup P_4,\text{bull},F_1)$-free graph is perfectly weight divisible. Suppose, to the contrary, that $G$ is a minimal nonperfectly weight divisible $(P_2\cup P_4,\text{bull},F_1)$-free graph. By the minimality of $G$, the graph $G$ must be connected. By Lemma~\ref{homogeneous},
 	\begin{equation}\label{eq3-1}
 		G~\text{has no homogeneous set}.
 	\end{equation}
 	
 	Moreover, we prove that, for every $x\in V(G)$,
 	\begin{equation}\label{eq3-2}
 		\mbox{$G[N(x)]$ is perfect and $G[M(x)]$ is imperfect.}
 	\end{equation}
 	
 	By (\ref{eq3-1}) and Lemma~\ref{bullfree}, either $G[N(x)]$ or $G[M(x)]$ is perfect. Suppose that $G[M(x)]$ is perfect. Then $G[\{x\}\cup M(x)]$ is perfect, since $x$ is anticomplete to $M(x)$. Moreover, for every positive integer weight function $w$ on $V(G)$, we have $\omega_w(G[N(x)])<\omega_w(G)$. Indeed, if $K$ is a maximum-weight clique of $G[N(x)]$, then $K\cup\{x\}$ is a clique of $G$, and hence $\omega_w(G)\geqslant\omega_w(G[N(x)])+w(x)>\omega_w(G[N(x)])$. Thus, $(\{x\}\cup M(x),N(x))$ is a partition witnessing that $G$ is perfectly weight divisible, a contradiction. Therefore, $G[M(x)]$ is imperfect, and consequently $G[N(x)]$ is perfect by Lemma~\ref{bullfree}. This proves (\ref{eq3-2}).
 	
 	First, we consider the case where $\omega(G)\leqslant2$. In this case, $\chi(G)\leqslant3$ by Corollary~\ref{co-1}. Consequently, $G$ is perfectly weight divisible by Lemma~\ref{3-coloring}, a contradiction. Therefore, $\omega(G)\geqslant3$. Let $v\in V(G)$ be a vertex contained in a maximum clique of $G$. By (\ref{eq3-2}), $G[N(v)]$ is perfect and $G[M(v)]$ is imperfect. We next prove the following claim.
 	
 	\begin{claim}\label{cl3-1}
 		$G[M(v)]$ contains a $5$-hole.
 	\end{claim}
 	
 	\pf
 	Suppose, to the contrary, that $G[M(v)]$ contains no $5$-hole. Since $G[M(v)]$ is imperfect by (\ref{eq3-2}), Theorem~\ref{perfect} implies that $G[M(v)]$ contains an odd hole or an odd antihole. Since every odd hole of length at least $9$ contains an induced $P_2\cup P_4$, $G[M(v)]$ contains either a $7$-hole or an odd antihole of length at least $7$.
 	
 	We first show that $G[M(v)]$ is $7$-hole-free. Suppose that $C'$ is a $7$-hole in $G[M(v)]$. By (\ref{eq3-2}), $G[N(x)]$ is perfect for every $x\in V(G)$. Hence, no vertex of $G$ is complete to $V(C')$, since otherwise its neighborhood would contain the $7$-hole $C'$ as an induced subgraph. By Lemma~\ref{c7}, $V(C')$ is a dominating set of $G$. However, $C'\subseteq M(v)$, and hence $v$ is anticomplete to $V(C')$, a contradiction. Therefore, $G[M(v)]$ is $7$-hole-free.
 	
 	It follows that $G[M(v)]$ contains an odd antihole $H$ with $V(H)=\{v_1,v_2,\ldots,v_k\}$, where $k\geqslant7$ is odd and $\overline{H}=v_1v_2\cdots v_kv_1$. Let $v'\in N(v)$. We first prove that
 	\begin{equation}\label{eq3-3}
 		|N(v')\cap V(H)|\geqslant2.
 	\end{equation}
 	
 	If $N(v')\cap V(H)=\emptyset$, then $\{v,v',v_1,v_3,v_k,v_2\}$ induces a $P_2\cup P_4$. If $|N(v')\cap V(H)|=1$, without loss of generality, let $N(v')\cap V(H)=\{v_1\}$. Then $\{v,v',v_3,v_5,v_2,v_4\}$ induces a $P_2\cup P_4$. Both are contradictions. This proves (\ref{eq3-3}).
 	
 	Next, we prove that
 	\begin{equation}\label{eq3-4}
 		N(v')\cap V(H)~\text{is a stable set}.
 	\end{equation}
 	
 	Suppose, to the contrary, that $N(v')\cap V(H)$ is not a stable set. Without loss of generality, we may suppose that $v_1,v_n\in N(v')\cap V(H)$ with $v_1v_n\in E(G)$, where $3\leqslant n\leqslant k-1$. Suppose that $v'$ is not complete to $\{v_1,v_2,\ldots,v_n\}$. Let $2\leqslant n'\leqslant n-1$ be the minimum integer such that $v'\not\sim v_{n'}$.
 	
 	If $n=3$, then $n'=2$. To avoid an induced bull on $\{v',v_1,v_3,v,v_4\}$, we have $v'\sim v_4$; but then $\{v',v_1,v_4,v,v_2\}$ induces a bull, a contradiction. Hence, $n\geqslant4$, and thus $v_n\sim v_2$ and $v_{n-1}\sim v_1$. Consequently, $n'\neq2$ to avoid an induced bull on $\{v',v_1,v_n,v,v_2\}$, and $n'\neq n-1$ to avoid an induced bull on $\{v',v_1,v_n,v,v_{n-1}\}$. Thus, $3\leqslant n'\leqslant n-2$, which implies that $v_{n'}\sim v_n$ and $v'\sim v_{n'-1}$ by the minimality of $n'$. But then $\{v',v_{n'-1},v_n,v,v_{n'}\}$ induces a bull, a contradiction. Therefore, $v'$ is complete to $\{v_1,v_2,\ldots,v_n\}$.
 	
 	Similarly, we deduce that $v'$ is complete to $\{v_1,v_k,v_{k-1},\ldots,v_n\}$. Hence, $v'$ is complete to $V(H)$. It follows that $G[N(v')]$ contains the odd antihole $H$ as an induced subgraph and is therefore imperfect, contradicting (\ref{eq3-2}). This proves (\ref{eq3-4}).
 	
 	Combining (\ref{eq3-3}) and (\ref{eq3-4}), and observing that every stable set of the odd antihole $H$ has size at most $2$, we have $|N_H(v')|=2$. Moreover, these two vertices are consecutive in $\overline{H}$. Thus, without loss of generality, we may assume that $N_H(v')=\{v_1,v_2\}$. But then $\{v,v',v_4,v_6,v_3,v_5\}$ induces a $P_2\cup P_4$, a contradiction. This proves Claim~\ref{cl3-1}. \qed
 	
 	By Claim~\ref{cl3-1}, let $C=v_1v_2v_3v_4v_5v_1$ be a $5$-hole in $G[M(v)]$. By Lemma~\ref{c5}, for every vertex $u\in N(v)$, either $N_C(u)=\{v_i,v_{i+2}\}$ for some $i\in[5]$, or $u$ is complete to $V(C)$. The latter is impossible by (\ref{eq3-2}), since otherwise $G[N(u)]$ would contain the $5$-hole $C$ as an induced subgraph and hence would be imperfect. Therefore,
 	\begin{equation}\label{ee-1}
 		\mbox{for every vertex $u\in N(v)$, $N_C(u)=\{v_i,v_{i+2}\}$ for some $i\in[5]$.}
 	\end{equation}
 	
 	From now on, the subscripts are taken modulo $5$. We prove the following claim.
 	
 	\begin{claim}\label{cl3-2}
 		Let $u,u'\in N(v)$ such that $u\sim u'$. Then $N_C(u)=N_C(u')$.
 	\end{claim}
 	
 	\pf
 	Suppose, to the contrary, that $N_C(u)\neq N_C(u')$. Without loss of generality, let $N_C(u)=\{v_1,v_3\}$ by (\ref{ee-1}). If $N_C(u')=\{v_2,v_4\}$, then $\{v,u,u',v_1,v_4\}$ induces a bull. If $N_C(u')=\{v_3,v_5\}$, then $\{u,u',v_3,v_2,v_5\}$ induces a bull. The remaining cases are symmetric. Hence, $N_C(u)=N_C(u')$. This proves Claim~\ref{cl3-2}. \qed
 	
 	Recall that $\omega(G)\geqslant3$ and $v$ belongs to a maximum clique of $G$. Thus, there exist two adjacent vertices in $N(v)$. Let $S$ be a component of $G[N(v)]$ that contains an edge. Since $S$ contains an edge and $v\notin V(S)$, we have $1<|V(S)|<|V(G)|$. We next prove the following claim.
 	
 	\begin{claim}\label{homo}
 		$V(S)$ is a homogeneous set of $G$.
 	\end{claim}
 	
 	\pf
 	Suppose that $V(S)$ is not a homogeneous set of $G$. By (\ref{ee-1}) and Claim~\ref{cl3-2}, and since $S$ is connected, we may assume by symmetry that $N(V(S))\cap V(C)=\{v_1,v_3\}$.
 	
 	Clearly, $\{v,v_1,v_3\}$ is complete to $V(S)$, and $(N(v)\setminus V(S))\cup\{v_2,v_4,v_5\}$ is anticomplete to $V(S)$. Since $V(S)$ is not a homogeneous set of $G$, there exists a vertex $z$ that has both a neighbor and a nonneighbor in $V(S)$. Clearly, $z\notin N(v)\cup\{v\}\cup V(C)$. Since $S$ is connected, there exists an edge $xy$ of $S$ such that $z\sim x$ and $z\not\sim y$.
 	
 	We next show that
 	\begin{equation}\label{e-111}
 		\mbox{$z$ is anticomplete to $\{v_1,v_3\}$.}
 	\end{equation}
 	
 	Suppose, to the contrary, that $z$ has a neighbor in $\{v_1,v_3\}$. By symmetry, we may assume that $z\sim v_1$. To avoid an induced bull on $\{x,z,v,v_1,v_2\}$ or $\{x,z,v,v_1,v_5\}$, we have that $z$ is complete to $\{v_2,v_5\}$. To avoid an induced bull on $\{z,y,v_1,v_4,v_5\}$, we have $z\sim v_4$. But now $\{z,y,v_1,v_2,v_4\}$ induces a bull, a contradiction. This proves (\ref{e-111}).
 	
 	To avoid an induced bull on $\{x,y,z,v_1,v_2\}$ or $\{x,y,z,v_3,v_4\}$, we have $z\sim v_2$ and $z\sim v_4$. Similarly, $z\sim v_5$. But now $\{z,v_1,v_3,v_4,v_5\}$ induces a bull, a contradiction. Therefore, $V(S)$ is a homogeneous set of $G$. This proves Claim~\ref{homo}. \qed
 	
 	\medskip
 	
 	By Claim~\ref{homo}, $V(S)$ is a homogeneous set of $G$, contradicting Lemma~\ref{homogeneous}. This completes the proof of Theorem~\ref{perfectdivisible}. \qed
 	\medskip
 	
 	
 	
 	{\small	\section*{Declarations}
 		\begin{itemize}
 			\item \noindent{\bf Acknowledgements} The authors would like to thank Dr. Ran Chen for bringing this problem to their attention and for his valuable suggestions.
 			\item \textbf{Author Contributions}\quad The authors contributed equally.
 			\item \textbf{Conflict of interest}\quad The authors declare no conflict of interest.
 			\item \textbf{Data availability statement}\quad This manuscript has no associated data.
 	\end{itemize}}

 \end{document}